%% file: arxiv.tex
\documentclass[preprint,12pt]{elsarticle}

\usepackage{amssymb}
\usepackage{amsmath}
\usepackage{color}
\usepackage{theorem}
\usepackage{subfigure}
 \usepackage{multirow}
 \usepackage{colortbl, xcolor}
 \definecolor{Gray}{gray}{0.9}

\usepackage[toc,page]{appendix}

\newtheorem{theorem}{Theorem}
\newtheorem{lemma}{Lemma}

\def\e1{1\!\!1}
\def\ex{\mathbf{x}}
\newcommand{\PP}{\mathbb{P}}
\newcommand{\R}{\mathbb{R}}

\newcommand{\E}{\mathbb{E}\,}

\def\eb{\textrm{\mathversion{bold}$\mathbf{\beta}$\mathversion{normal}}}  

\def\eth{\textrm{\mathversion{bold}$\mathbf{\theta}$\mathversion{normal}}}    

\journal{Econometrics and Statistics, Part B}

\begin{document}

\begin{frontmatter}




\title{Quantile LASSO with changepoints in panel data models applied to option pricing}


\author{Mat\'u\v{s} Maciak}

\ead{maciak@karlin.mff.cuni.cz}


\address{Department of Probability and Mathematical Statistics, Faculty of Mathematics and Physics, Charles University, Sokolovsk\'a 83, Prague, 186 75, Czech Republic}

\begin{abstract}
Panel data are modern statistical tools which are commonly used in all kinds of econometric problems under various regularity assumptions. The panel data models with changepoints are introduced together with atomic pursuit methods and they are applied to estimate the underlying option price function. Robust estimates and complex insight into the data are both achieved by adopting the quantile LASSO approach. The final model is produced in a fully data-driven manner in just one single modeling step. In addition, the arbitrage-free scenarios are obtained by introducing a set of well defined linear constraints. The final estimate is, under some reasonable assumptions, consistent with respect to the model estimation and the changepoint detection performance. The finite sample properties are investigated in a simulation study and proposed methodology is applied for the Apple call option pricing problem.
\end{abstract}


\end{frontmatter}


\section{Introduction}
\label{introduction}
The panel data and changepoints are frequently discussed and hot topics in theoretical and empirical econometrics. On the other hand, the option pricing problem and the corresponding implied volatility surface estimation, which both represent one of the most fundamental problem in financial mathematics and derivatives trading (see \cite{britten}), are mostly based on various approaches derived from the well-known Black-Scholes model \cite{BlackScholes73}. This model is still popular among practitioners that it is considered to be unrealistic from the theoretical point of view. In this paper we combine these two areas and we propose the option pricing technique based on the idea of the panel data models and the sparse estimation via atomic pursuit methods. A similar approach was recently studied in \cite{qian_su} but the authors only considered the standard quadratic loss function. In our approach we use the quantile check function which is, in general, more robust and it offers a more complex insight into the data as it can estimate an arbitrary conditional quantile rather than just the conditional mean.  

From the theoretical point of view, our method is motivated by the concept of  a regularized changepoint detection proposed in \cite{Harchaoui.Levy.10} and further elaborated for the conditional quantile estimation in \cite{Ciuperca-Maciak-17} or \cite{Ciuperca-Maciak-19}.  The main advantage of our approach relies on four pivots: Firstly, the method is robust respect to the noise in the price observations caused by  various market effects (for instance, bid-ask spreads, discrete ticks in price, non-synchronous trading); Second, the optimization problem is convex and  the existence of the optimal solution is, therefore, guaranteed; Third, the fully automated estimation process is used  with no need for any nuisance parameters to be pre-determined; Finally, the overall simplicity allows for the arbitrage-free scenarios which are obtained in a straightforward way by  a set of some well imposed linear constraints.

 From the practical point of view, our method is motivated by various semiparametric and nonparametric option pricing approaches used in econometrics (see, for instance, \cite{Benko07}, \cite{Fengler05}, or \cite{Kahale04})  where the option price function or the corresponding implied volatility function are usually obtained in terms of some constrained minimization problem. Alternatively, the option price function based can be also directly used to interpolate the implied volatility surface (see \cite{Fengler05} or  \cite{hull}).

The rest of the paper is organized as follows: the quantile fussed LASSO model is briefly introduced in Section \ref{fused_lasso}. For illustration, the model is applied in Section \ref{fixed_maturity} to estimate  the option price function for a fixed time and a given maturity. The model is later generalized for the panel data structure
in Section \ref{all_maturities} and this model is applied for the Apple call option pricing problem in Section \ref{application}. A small simulation study is also presented in Section \ref{application}. 

\section{Quantile fused LASSO}
\label{fused_lasso}
Consider a standard linear model, however, with the parameters which can change along the available observations $i \in \{1, \dots, n\}$, such that
\begin{equation}
\label{eq1}
Y_i =\ex_i^\top \eb_i+\varepsilon_i, \qquad i=1, \cdots , n,
\end{equation}
where  $\eb_i \in \R^p$ is a $p$-dimensional parameter (the dimension does not depend on $n \in \mathbb{N}$) and $\ex_i=(x_{i1}, x_{i2}, \cdots , x_{ip})^\top$ is the subject's specific vector of explanatory variables. The error terms $\{\varepsilon_i\}_{i = 1}^n$  are supposed to be independent, centered, and identically distributed with some generally unknown distribution function $F$. In addition, we assume that  there is a specific sparsity structure  in $\eb_i$'s, such that $\eb_i = \eb_{i - 1}$, for most of the indexes $i \in \{2, \dots, n\}$, but some few exceptions---changepoints.
Such model can be seen as a straightforward extension of a simple piece-wise constant model discussed in \cite{Harchaoui.Levy.10} or, from the econometrics perspective, a generalization of a more common trend model proposed in \cite{maciak}. The same model as in \eqref{eq1}, however, for the dependent time series data, is also considered in \cite{qian_su2}. 

The model in \eqref{eq1} is assumed to have $K^* \in \mathbb{N}$ changepoints in total,  located at $t^*_1 < \cdots < t^*_{K^*} \in \{1, \dots, n\}$, such that 
\begin{equation}
\label{eq2}
\eb_i=\eb_{t_k^*}, \qquad \forall i=t^*_k, t_k^*+1, \cdots , t^*_{k+1}-1, \qquad k=0,1, \cdots , K^*, 
\end{equation}
with $t^*_0=1$, $t^*_{K^*+1}=n$, and $\eb_n=\eb_{t^*_{K^*+1}}$.  In general, the number of true changepoints $K^* \in \mathbb{N}$ and their locations $t^*_1, \cdots , t^*_{K^*}$ are all unknown. The true values of $\eb_i$ are denoted by $\eb^*_i$ and $K^* \equiv Card\{ i \in \{2, \cdots , n \}; \; \eb^*_i \neq \eb^*_{i-1}  \}$. 
The idea of the estimation method is to recover the unknown changepoint locations and to estimate the underlying model phases---the vector parameters which are associated with the conditional quantiles of interest. For this purpose, the following optimization problem is formulated

\begin{equation}
\label{minimization1}
\widehat{\eb^n} =
\begin{array}{c}
~\\[-0.4cm]
\textrm{Argmin}\\[-0.2cm]
{\scriptsize \eb_i \in \mathbb{R}^{p};~ i = 1, \dots, n}
\end{array} \quad \sum^n_{i=1} \rho_\tau(Y_i- \ex_i^\top \eb_i)+n\lambda_n \sum^n_{i=2} \|\eb_i-\eb_{i-1}\|_2,
\end{equation}
where, for simplicity, $\widehat{\eb^n}=(\widehat{\eb}_1^\top, \cdots, \widehat{\eb}_n^\top)^\top \in \mathbb{R}^{np}$, $\rho_{\tau}(u) = u (\tau - \mathbb{I}_{\{u < 0\}})$, for $\tau \in (0,1)$, is the standard check function used for the quantile regression, $\|\cdot\|_2$ stands for the classical $L_2$ norm, and $\lambda_n >0 $ is the tuning parameter which controls for the overall number of changepoints (the sparsity level) occurring in the final model: for $\lambda_{n} \to 0$ there will be $\widehat{\eb}_i \neq \widehat{\eb}_{i - 1}$ for each $i \in \{2, \dots, n\}$, while for $\lambda_n \to \infty$ 
 no changepoints are expected to occur in the final model and, thus, $\widehat{\eb}_i = \widehat{\eb}_{i - 1}$ for all $i \in \{2, \dots, n\}$. The corresponding estimators for the changepoint  locations are the observations $i \in \{2, \cdots , n \}$, where $\widehat{ \eb}_i \neq \widehat{\eb}_{i-1} $. Let us, therefore,  define the set 
\begin{equation}
\label{hAn}
\widehat{ {\cal A}}_n \equiv \{ i \in \{2, \cdots , n \}; \; \widehat{ \eb}_i \neq \widehat{ \eb}_{i-1}  \} = \{ \hat{  t}_1 < \cdots < \hat{  t}_{|\widehat{ \cal A}_n|} \},
\end{equation}
and let $|\widehat{ \cal A}_n|$ be the cardinality of $\widehat{ {\cal A}}_n$. For each $k = 0, \dots, |\widehat{ {\cal A}}_n|$ we can also define the $(k + 1)$-st model phase (observations indexed by the set $\{\hat{t}_k, \dots, \hat{t}_{k + 1} - 1\}$, where $\hat{t}_0 =1$ and $\hat{t}_{\widehat{ {\cal A}}_n + 1} = n$), with the corresponding vector of estimated parameters $\widehat{\eb}_{\hat{t}_{k}}$.  The minimization problem formulated in \eqref{eq1} is convex and it can be effectively solved using some standard optimization toolboxes (see, for instance, \cite{Huang}).    The theoretical properties are studied in detail in \cite{Ciuperca-Maciak-19}. Under some reasonable assumptions, the method achieves consistency in terms of the changepoint detection and, also, in terms of the parameter estimation. Nevertheless, the regularization parameter in the LASSO problems should be, in general, chosen differently when aiming at the changepoint recovery or the underlying model estimation: for the former one,  larger values are preferred to avoid the overestimation issue and false changepoint detection. On the other hand, for the estimation purposes, slightly smaller values of $\lambda_n > 0 $ are needed in order to limit the shrinkage effect and to improve the estimation bias performance. The value of $\lambda_n > 0$ which satisfies the set of assumptions used in \cite{Ciuperca-Maciak-19} is, for instance, $\lambda_n = (1/n) \cdot (\log n)^{5/2}$. 

The role of the regularization parameter is crucial but its importance can be suppressed by using some alternative source of the regularization. This is, for instance, the case for the option pricing problem where  the final model must satisfy some shape constraints to comply with the financial theory on the arbitrage-free markets. Prescribing the convex and non-increasing property for the final estimate serves as an alternative regularization and thus, the choice of $\lambda_n > 0$ becomes rather inferior. This is also demonstrated in the next section.

\section{Option price function for a fixed maturity}
\label{fixed_maturity}
Let us start with a simple situation where the call option prices are observed for some specific maturity at some fixed time. The data can be represented as $\{(Y_i, x_{i});~i = 1, \dots, n\}$, where $Y_{i}$ stands for the option intrinsic  value at the strike $x_{i}$. In total, there are $n \in \mathbb{N}$ observations  and the aim is to estimate the option price function, which must be non-increasing and convex. The quantile fused LASSO presented in Section \ref{fused_lasso} is adopted to construct the estimate, however, a proper modification is needed to meet all desired qualities of the final model. Firstly, the quantile level of $\tau = 0.5$ is used and, thus, the conditional median will be obtained as the solution of \eqref{eq1}. However, if the density of the error terms is symmetric, the conditional median will, under some moment conditions,  coincide with the conditional mean. Moreover, the estimate will be robust with respect to possible outliers which is a convenient property in the derivatives trading. Second, the monotonic and isotonic properties are not automatically guaranteed by minimizing \eqref{eq1}, therefore, some additional constraints are needed. The desired qualities can be, however, obtained in a~straightforward way by introducing a set of well defined linear constraints: the resulting minimization problem is still convex and  the same algorithm can be again used to obtain the solution.

Let us assume that the strikes $\{x_{i}\}_{i = 1}^{n}$ are all from some compact domain, denoted as $\mathcal{D}$. Let $\{\varphi_j(x);~j = 1, \dots, p\}$ be some functional basis constructed on $\mathcal{D}$  and let $\boldsymbol{x}_{i} = (\varphi_1(x_i), \dots, \varphi_p(x_i))^\top$. 
For instance, one can simply define $\boldsymbol{x}_{i} = (1, x_i, x_i^2)^\top$, for $i = 1, \dots, n$, where $p = 3$, to obtain a standard quadratic fit. The subject's specific parameters, $\boldsymbol{\beta}_{i} = (\beta_{i 1}, \beta_{i 2}, \beta_{i 3})^\top \in \mathbb{R}^3$, induce o lot of flexibility in the model, thus, the  final model would be too haphazard if no further restrictions on the parameter vectors were imposed.

For the subject's specific vectors $\boldsymbol{x}_i$, for $i = 1, \dots, n$, we can define the model matrix  
\[
\mathbb{X} = \left[  
\begin{array}{cccc}
\ex_1^\top & \textbf{0} & \cdots & \textbf{0 } \\
\textbf{0} & \ex_2^\top &  \cdots & \textbf{0 } \\
\vdots& \vdots & \ddots & \vdots \\
\textbf{0 } & \textbf{0 } & \ldots & \ex_n^\top
\end{array}
\right]
\]
and the overall model can be also expressed as 
\begin{align}
\label{eq3}
\boldsymbol{Y} = \mathbb{X} \boldsymbol{\beta}^{n} + \boldsymbol{\varepsilon},
\end{align}
where $\boldsymbol{Y} = (Y_1, \dots, Y_n)^\top$, $\boldsymbol{\varepsilon} = (\varepsilon_1, \dots, \varepsilon_n)^\top$, and $\boldsymbol{\beta}^{n} = (\boldsymbol{\beta}_1^\top, \dots, \boldsymbol{\beta}_n^\top)^\top \in \mathbb{R}^{np}$.
For the model in \eqref{eq3} we can directly use the minimization formulation in \eqref{minimization1} but the solution is, in general, not smooth and, moreover, the qualitative properties known from the financial theory on the arbitrage-free markets  are not automatically guaranteed. In particular, the option price function is supposed to be non-increasing and convex with respect to the strikes. These properties can be explicitly enforced by minimizing \eqref{minimization1} with respect to some well defined linear constraints. Specifically, the final solution (the estimated option price function)  will be non-increasing if the estimated vector of parameters $\widehat{\boldsymbol{\beta}^n} \in \mathbb{R}^{np}$ obeys
\begin{equation}
\label{eq4}
\boldsymbol{D} \widehat{\boldsymbol{\beta}^n} \leq \boldsymbol{0},
\end{equation}
which holds element-wise for
$$
\boldsymbol{D} = \left[  
\begin{array}{cccc}
\widetilde{\ex}_1^\top & \textbf{0} & \cdots & \textbf{0 } \\
\textbf{0} & \widetilde{\ex}_2^\top &  \cdots & \textbf{0 } \\
\vdots& \vdots & \ddots & \vdots \\
\textbf{0 } & \textbf{0 } & \ldots & \widetilde{\ex}_n^\top
\end{array}
\right],
$$
where $\widetilde{\ex}_i = (\varphi_1'(x_{i}), \dots, \varphi_p'(x_i))^\top$ is the vector of the first derivatives of the functional basis $\{\varphi_j(x);~j = 1, \dots, p\}$ evaluated at the strike $x_i \in \mathcal{D}$, for $i = 1, \dots, n$. For some reasonable degrees of the functional basis the condition in \eqref{eq4} ensures that the final model is, indeed, non-increasing in the strikes. Analogously, the convexity of the final model (the estimated option price function) is enforced if the vector of the estimated parameters $\widehat{\boldsymbol{\beta}^n} \in \mathbb{R}^{np}$ satisfies
\begin{equation}
\label{eq5}
\boldsymbol{C} \widehat{\boldsymbol{\beta}^n} \geq \boldsymbol{0},
\end{equation}
which again holds element-wise for 
$$
\boldsymbol{C} = \left[  
\begin{array}{cccc}
\widetilde{\widetilde{\ex}}_1^\top & \textbf{0} & \cdots & \textbf{0 } \\
\textbf{0} & \widetilde{\widetilde{\ex}}_2^\top &  \cdots & \textbf{0 } \\
\vdots& \vdots & \ddots & \vdots \\
\textbf{0 } & \textbf{0 } & \ldots & \widetilde{\widetilde{\ex}}_n^\top
\end{array}
\right],
$$
where $\widetilde{\widetilde{\ex}}_i = (\varphi_1''(x_{i}), \dots, \varphi_p''(x_i))^\top$ denotes now the vector of the second derivatives of the functional basis $\{\varphi_j(x);~j = 1, \dots, p\}$ evaluated again at the strike $x_i \in \mathcal{D}$, for $i = 1, \dots, n$. The minimization problem in \eqref{minimization1} together with the shape constraints---the non-decreasing property enforced by \eqref{eq4}, and the convexity imposed by \eqref{eq5} is convex  and it can be effectively solved using the standard optimization toolboxes. 

\begin{figure}
	\centering
	\subfigure[\scriptsize No shape restrictions]{\label{fig1a}\includegraphics[width=0.32\textwidth, height = 0.24\textwidth]{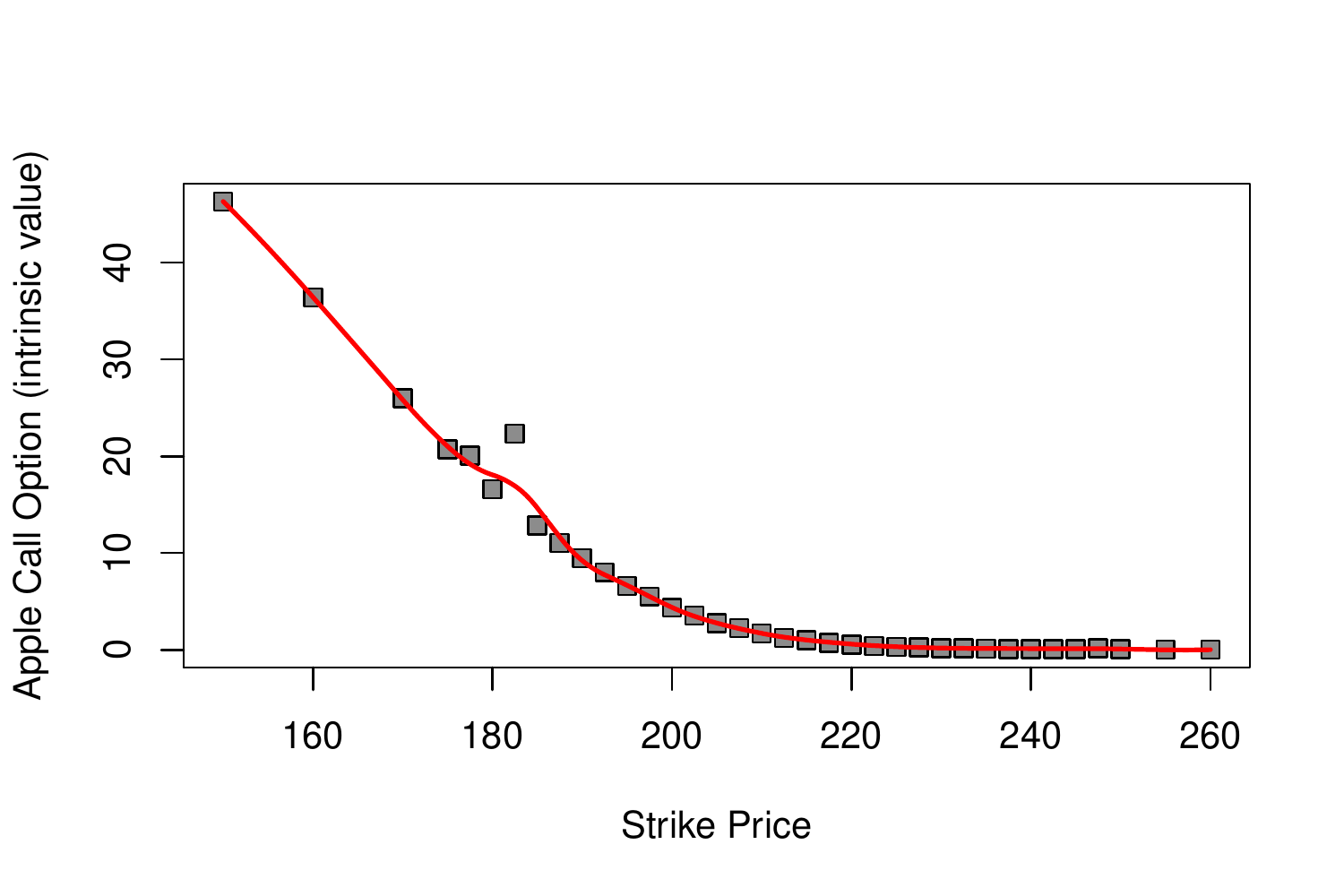}}
	\subfigure[\scriptsize Non-increasing property]{\label{fig1b}\includegraphics[width=0.32\textwidth, height = 0.24\textwidth]{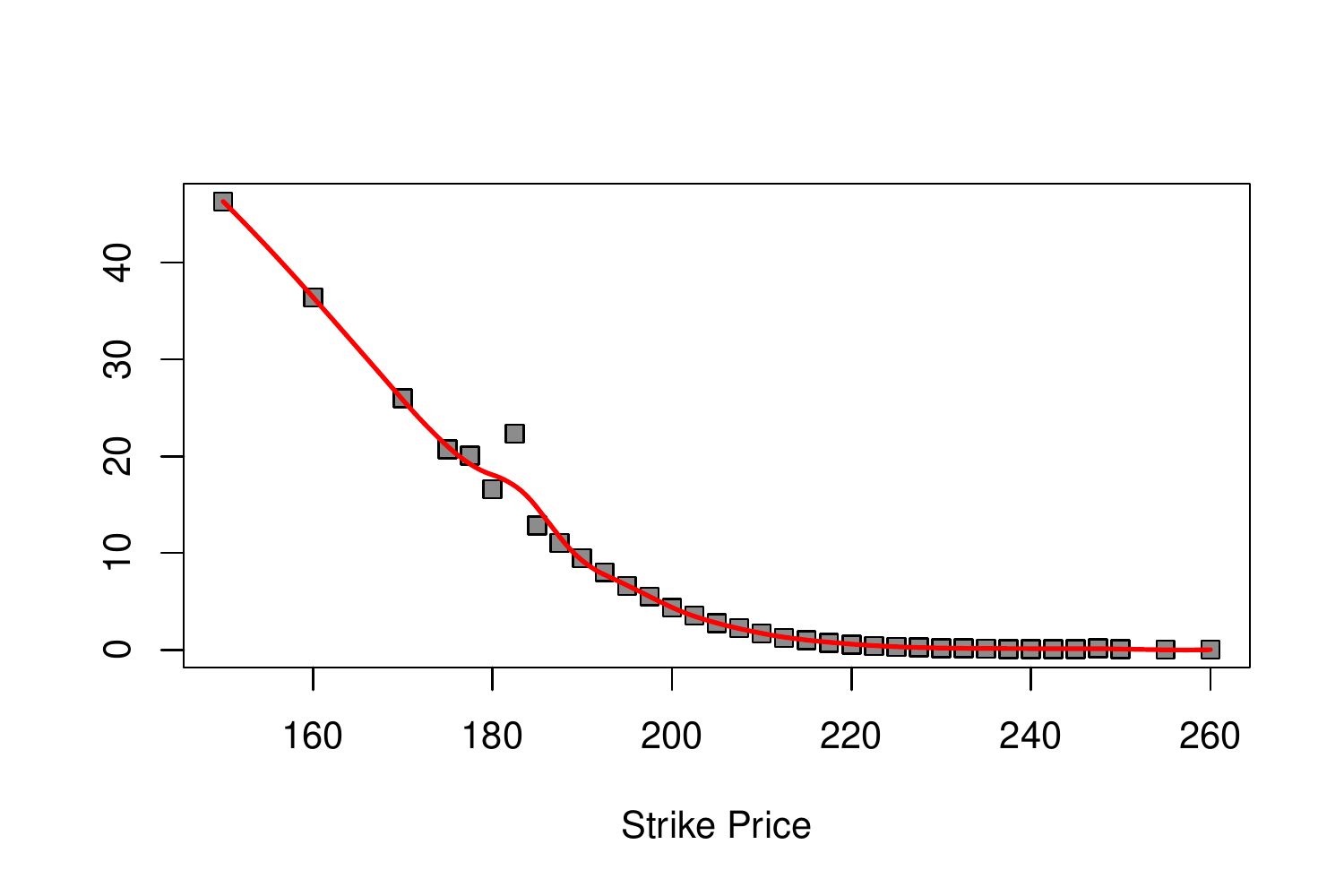}}
	\subfigure[\scriptsize Non-increasing and convex]{\label{fig1c}\includegraphics[width=0.32\textwidth, height = 0.24\textwidth]{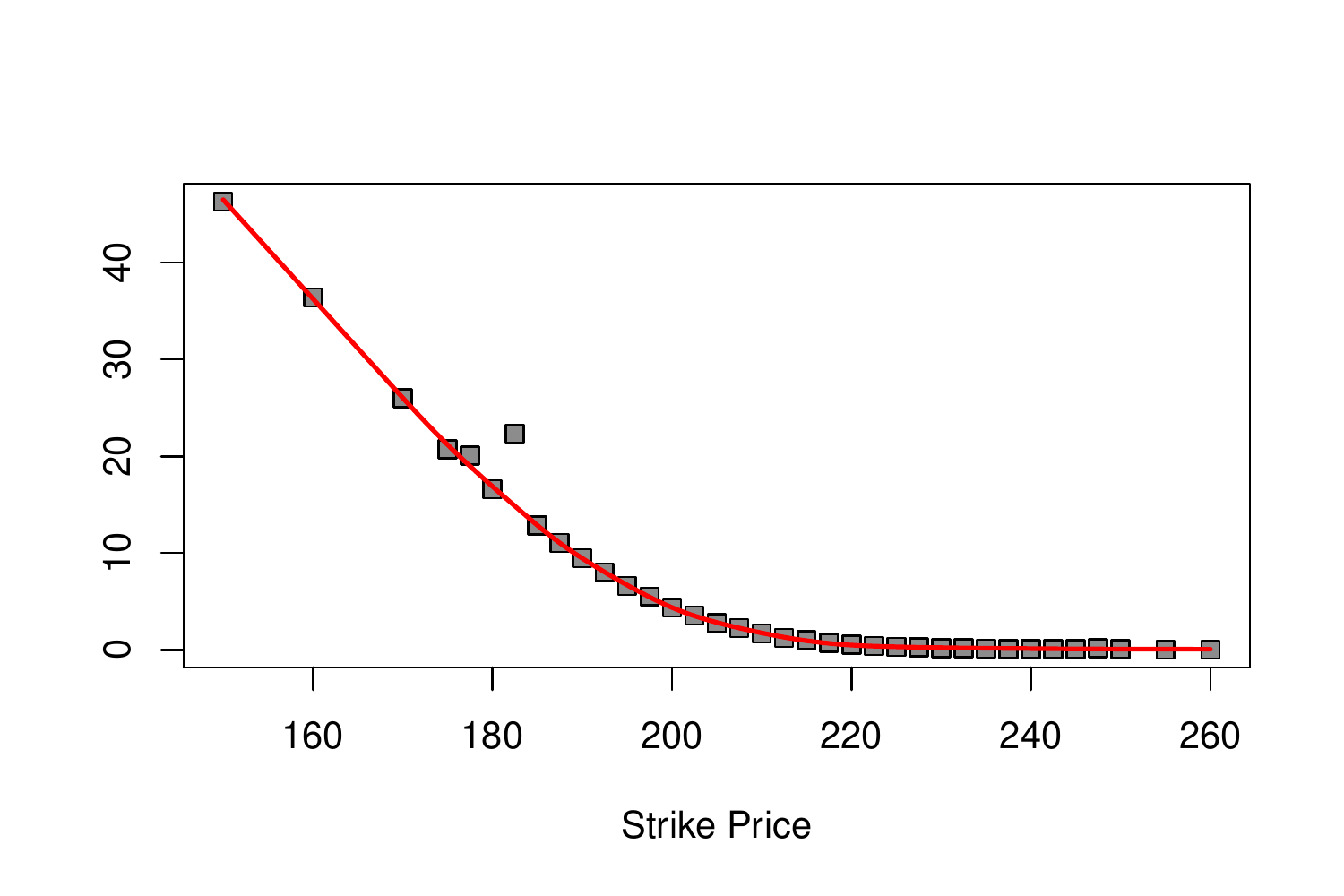}}
	
	\caption{\footnotesize The illustration of the Quantile LASSO performance when applied to the Apple Call Options with the maturity 32 days. In the left figure, the price function is estimated without enforcing any shape constraints. In the middle figure, the non-increasing property is used, however, without requiring convexity. Finally, the right figure shows the option price function while enforcing the non-increasing property and also convexity with respect to the strikes. The value of the regularization parameter   is $\lambda_n = n^{-1} (\log n)^{5/2}$ which is an example of a sequence which satisfies the assumptions in \cite{Ciuperca-Maciak-19}.}
	\label{fig1}
\end{figure}

For illustration, the proposed estimation approach is applied for the Apple Call Options with the expiration of 32 days (see Figure \ref{fig1}). The smoothness property is automatically achieved by using the third degree spline basis $\{\varphi_j(x);~j = 1, \dots, p\}$.  The non-increasing property seems to be automatically obtained by the nature of the data (see Figure \ref{fig1a})) without explicitly enforcing it. On the other hand, the convexity is more challenging and it is not directly guaranteed by the data (see Figure \ref{fig1a} or \ref{fig1b}) and it must be enforced (see Figure \ref{fig1c}) by using the shape constraints in \eqref{eq5}.

In practical applications, the estimated price function usually  changes over time as the time progresses towards the option's maturity---the expiry date. Because the expiry date is fixed, the follow-up period is limited and it can be relatively short. In the next chapter we introduce a modification of the quantile LASSO model which can be used for the panel data structure with $n \in \mathbb{N}$ panels observed over some relatively short time period $[0,T]$, for some fixed $T > 0$.

\section{Panel data model for time dependent maturity}
\label{all_maturities}
Let us now assume the data $\{(Y_{t i}, x_{i t});~t = 1, \dots, T;~i = 1, \dots, n\}$ where $Y_{i t}$ represents the option intrinsic value at some specific time $t \in \{1, \dots, T\}$ and the strike $x_{i t} \in \mathcal{D}$,  for $i = 1, \dots, n$. As far as the strikes are common over time, we can also assume that $x_{i t} \equiv x_i$, for all $i \in \{1, \dots, n\}$. 
Thus, for each available strike $x_i \in \mathcal{D}$ we have a  strike specific panel of the strike specific intrinsic values observed over time $t \in \{1, \dots, T\}$, for some fixed $T \in \mathbb{N}$. The underlying panel data model can be expressed as
\begin{equation}
\label{eq_panelmodel}
Y_{t i} = \boldsymbol{x}_i^\top\boldsymbol{\beta}_t + \varepsilon_{t i}, \quad \textrm{for $t = 1, \dots, T$ and $i = 1, \dots, n$,}
\end{equation} 
where again $\boldsymbol{x}_i = (\varphi_1(x_i), \dots, \varphi_p(x_i))^\top$ is some functional basis on $\mathcal{D}$ evaluated at the given strike, $\boldsymbol{\beta}_t = (\beta_{t 1}, \dots, \beta_{t p})^\top \in \mathbb{R}^p$ is the vector of unknown parameters which can change over time $t \in \{1, \dots, T\}$ and the error vectors $\boldsymbol{\varepsilon}_i = [\varepsilon_{1 i}, \dots, \varepsilon_{T i}]$ are independently distributed over $i \in \{1, \dots, n\}$.

In order to estimate all panels simultaneously, such that the final model will obey the shape restrictions required for the arbitrage-free market, the following minimization problem is considered
 
\begin{align}
\label{surface-estimation}
\begin{array}{c}
~\\[-0.4cm]
\textrm{Minimize}\\[-0.2cm]
{\scriptsize \eb_{t} \in \mathbb{R}^{p}}
\end{array} & \hskip0.2cm \sum_{t = 1}^{T} \sum_{i = 1}^{n} \rho_\tau \Big(Y_{t i} - \boldsymbol{x}_{i}^\top \boldsymbol{\beta}_{t} \Big) + n \lambda_{N} \sum_{t = 2}^{T}  \|\eb_{t}-\eb_{t - 1}\|_2
\end{align}
with respect to
\begin{align}
\label{constraints}
\left.
\begin{array}{lrc}
\boldsymbol{D}\boldsymbol{\beta}_{t} \leq \boldsymbol{0},~t = 1, \dots, T; & \textrm{\it (non-increasing in the strike)} & (C1)\\
 \boldsymbol{C}\boldsymbol{\beta}_{t} \geq \boldsymbol{0},~t = 1, \dots, T; & \textrm{\it (convexity in the strike)}& (C2) 
\end{array}\right\}
\end{align}
where the matrices $\boldsymbol{D}$ and $\boldsymbol{C}$ are defined analogously as in \eqref{eq4} and \eqref{eq5} respectively. The overall vector of the estimated parameters $\widehat{\boldsymbol{\beta}^n} = (\widehat{\boldsymbol{\beta}}_{1}^\top, \dots, \widehat{\boldsymbol{\beta}}_T^\top)^\top \in \mathbb{R}^{T \times p}$ represents the set of all panels while $\widehat{\boldsymbol{\beta}}_t \in \mathbb{R}^{p}$ is only associated with the estimated of the option price function for the specific time $t \in \{1, \dots, T\}$. The price function is obviously allowed to change over time to reflect possible changes at the market but due to the fused penalty term in \eqref{surface-estimation} the model only allows for some of these changes to occur in the final model. The following sparsity structure is assumed: for situations where $\widehat{\boldsymbol{\beta}}_t \neq \widehat{\boldsymbol{\beta}}_{t - 1}$ the option price function changes from time $(t -1)$ to time $t$ to adapt for the situation at the market, otherwise, the option price function remains the same. The regularization parameter $\lambda_n > 0$ controls the amount of such changes in the model and the shape constraints in \eqref{constraints} are responsible for the additional source of the regularization by enforcing the non-increasing and convex properties of the estimated option price function at each time point $t \in \{1, \dots, T\}$. The minimization problem in \eqref{surface-estimation} together with the set of the linear constraints in \eqref{constraints} is again a convex minimization problem and the optimal solution can be obtained by the standard optimization software. The Karush-Kuhn-Tucker (KKT) optimality conditions  can be easily derived and they are formulated by the following lemma. 

\begin{lemma}
	\label{lemma1}
	~\\[-0.8cm]
\begin{enumerate}[(a)]
	\item For any $l \in \{1, \dots, |\widehat{\mathcal{A}}_n|\}$, $n \in \mathbb{N}$, and $\lambda_n > 0$ the following holds with probability one: 
	$$\displaystyle{\tau (T - \hat{t}_l + 1) \sum_{i = 1}^n \ex_i -  \sum^n_{i = 1} \sum_{k = \hat{t}_l}^T \ex_i \e1_{\{Y_{i k} \leq \ex'_i \widehat{ \eb}_{k}\}}=n \lambda_n \frac{\widehat{ \eth}_{\hat{ t}_l}}{\| \widehat{\eth}_{\hat{t}_l}\|_2}};$$
	
	\item For any $t = \{1, \dots T\}$,  $n \in \mathbb{N}$, and $\lambda_n >0$, the following holds with probability one:
	 $$\displaystyle{\left\|\tau (T - t + 1)\sum^n_{i=1} \ex_i- \sum^n_{i=1} \sum_{k = t}^T \ex_i \e1_{\{Y_{i k} \leq \ex'_i \widehat{ \eb}_{k }\}} \right\|_2 \leq n \lambda_n}.$$
\end{enumerate}
\end{lemma}

The proof of Lemma \ref{lemma1} is straightforward: Using a reparametrization of the form $\eth_t = \eb_t - \eb_{t - 1}$, for $t = 2, \dots, T$, one just needs to realize that $\eb_t = \sum_{k = 1}^t \eth_k$, where $\eth_1 \equiv \eb_1$. The rest already follows directly from the definition of the quantile check function $\rho_{\tau}$. Using a similar idea as in the proofs in \cite{Ciuperca-Maciak-19} the KKT conditions above can be used to derive some theoretical properties. Let us just recall, that the follow-up period $T \in \mathbb{N}$ is assumed to be fixed and relatively short (compared to the number of panels, $n \in \mathbb{N}$, which are allowed to tend to infinity). Thus, the theoretical results are only derived for the situation where the number of panels increases. 
In the following we focus on some consistency results when estimating the true time-specific parameter $\eb_t$, for $t = 1, \dots, T$. Let us firstly state some necessary conditions which are required for the main result to hold.

\vskip0.4cm
\noindent\textbf{Assumptions:}
\begin{itemize}	
	\item[\textbf{(A1)}] The errors $\boldsymbol{\varepsilon}_i = [\varepsilon_{1 i}, \dots, \varepsilon_{T i}]$ are independent copies  of some strictly stationary sequence $\boldsymbol{\varepsilon} = [\varepsilon_1, \dots, \varepsilon_T]$ with the continuous marginal distribution functions $F_{\varepsilon_t}(x)$ and $F_{(\varepsilon_{t}, \varepsilon_{t + k})}(x, y)$, for $x, y \in \R$, $t, \in \{1, \dots, T\}$, and $k \geq 1$. Moreover, $F_{\varepsilon_t}(0) = \PP[\varepsilon_{t} <0]=\tau$, for $\tau \in (0,1)$. The corresponding density functions $f(\cdot)$ and $f(\cdot,\cdot)$ are  bounded and strictly positive in the neighborhood of zero;
	
	\item[\textbf{(A2)}] There exist two constants $c, C \in \mathbb{R}$ such that 
	$$0 < c \leq \mu_{min}(\E[\mathbb{X}_n]) \leq \mu_{max} (\E[\mathbb{X}_n]) \leq C < \infty,$$
	where $\mu_{min}$ and $\mu_{max}$ stand for the minimum and maximum eigenvalue of the matrix in the argument and $\mathbb{X}_n = \frac{1}{n} \sum_{i = 1}^n \boldsymbol{x}_i\boldsymbol{x}_i^\top$. Moreover, $\max_{1 \leq i \leq n} \|\boldsymbol{x}_{i} \|_\infty < \infty$.
	
	\item[\textbf{(A3)}] There are two deterministic positive sequences  $(\lambda_n)$ and $(\delta_n)$, such that $\lambda_n \to 0$,  $\delta_n \to 0$,  $ n^{1/2} \delta_n \to \infty$, and $\lambda_n / \delta_n \to 0$ as $n \to \infty$.
\end{itemize}

Let us recall that in similar models (see, for instance, \cite{Ciuperca-Maciak-19}, \cite{Ciuperca-Maciak-17}, or \cite{qian_su}) there is an additional assumption which requires that the span between two consecutive changepoints increases. Analogously, the overall number of chances in the model is usually considered to be fixed. However, as far as the follow-up period $T \in \mathbb{N}$ is assumed to be fixed these two assumptions are irrelevant for our situation. The main consistency results is given by the next theorem. 

\begin{theorem}
	\label{thm1}
	Let the assumptions in (A1)--(A3) be all satisfied. Then, for any $t = 1, \dots, T$, it holds that $$
	\|\widehat{\eb}_t - \eb_t^*\|_1 = O_{P}\left(\sqrt{\frac{\log n}{n}}\right),
	$$
	where $\widehat{\eb}_t \in \mathbb{R}^p$ denotes the vector of estimated parameters obtained by minimizing \eqref{surface-estimation} and $\eb_t^*$ is the corresponding vector of the true values. 
\end{theorem}

The theorem above shows the consistency achieved by the estimation procedure defined by \eqref{surface-estimation}. An example of two sequences, which satisfy Assumption (A3) are, for instance, $\lambda_n = n^{-1} \cdot (\log n)^{1/2}$ and $\delta_n = (n^{-1} \log n)^{1/2}$. The sketch of the proof is postponed to the appendix section.

Finally, let briefly discuss some options for the selection of the functional basis constructed on the strike domain $\mathcal{D}$. This basis implicitly occurs in \eqref{surface-estimation} via the term 
\begin{equation}
\label{basis}
\boldsymbol{x}_{i}^\top \boldsymbol{\beta}_{t} = \sum_{\ell = 1}^p \varphi_{\ell}(x_{i}) \beta_{t \ell},
\end{equation}
where $\boldsymbol{\beta}_{t} = (\beta_{t 1}, \dots, \beta_{t p})^\top \in \mathbb{R}^p$ is the time specific vector of unknown parameters which defines the option price function at the time $t \in \{1, \dots,T\}$. 
Natural candidates for the basis selection are, for instance, polynomials $\{\varphi_\ell(x) = x^\ell;~ \ell = 0, \dots, d\}$ of some degree $d \in \mathbb{N}$ (for instance, $d = 2$ or $d = 3$). Such basis, however, does not automatically enforce the overall continuity and smoothness of the final estimate. Therefore, truncated splines defined with respect to some set of inner knots in $\mathcal{D}$ are more appropriate. They can be easily defined in a way that the continuity and smoothness conditions are both achieved automatically and one just needs to minimize \eqref{surface-estimation} with respect to the linear constraints (C1)--(C2) in \eqref{constraints}. The selection of the inner knot points is not crucial either. The overestimation issue is controlled by the shape constraints in \eqref{constraints} which regularize the final fit even if the number of knots used for the spline basis is very large.

\section{Application to Apple option pricing}
\label{application}
\begin{figure}
	\centering
	\includegraphics[width=0.98\textwidth, height = 0.4\textwidth]{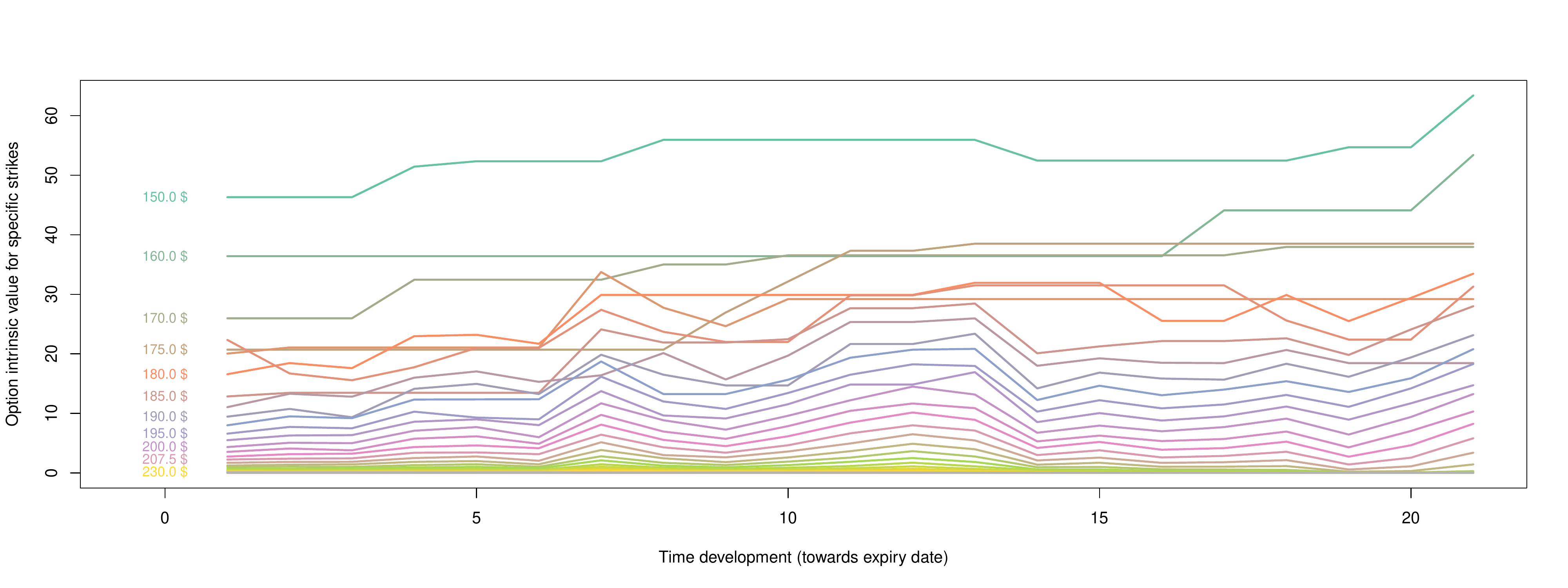}
	
	\caption{\footnotesize The daily development of the Apple Call option intrinsic values (from August 2019) with the initial maturity of 32 days given for 36 different strikes (panels). The intrinsic values are observed for 21 consecutive days (except weekends and holidays). The last observed values correspond with the option intrinsic values at the day of the maturity (September 6th, 2019). Given the financial theory for the call options, the intrinsic values are high for small strikes and they are close to zero for high strikes.}
	\label{fig1_alt}
\end{figure}

\begin{figure}
	\hskip-0.8cm\includegraphics[width=14cm, height = 6.5cm]{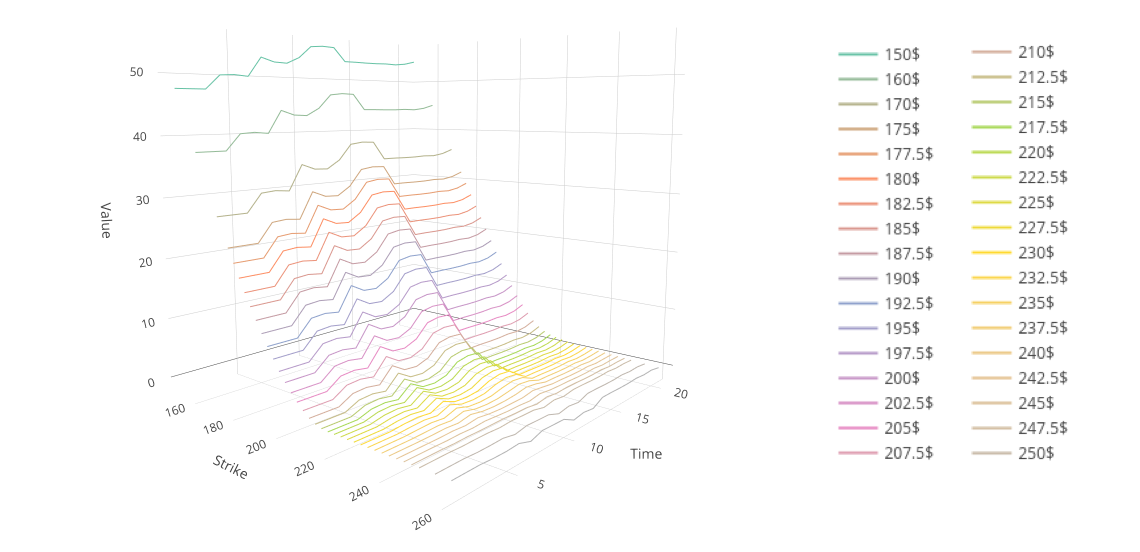}
	\caption{\footnotesize The estimated panels for the Apple Call Option price function based on 36 available strikes and the follow up period of 21 days (options with the initial maturity of 32 days observed daily except weekends and holidays). The estimated panels satisfy all qualitative properties required by the financial theory for the arbitrage-free market: the price function is always non-increasing and convex in the strikes at any given time point. It is also obvious that the price function slightly changes in the middle of the follow-up period. This change is allowed by the time-specific vectors of parameters used in \eqref{surface-estimation}. The corresponding time segments are $[1,3]$, $[4,6]$, $[7,10]$, $[11,13]$, and $[14,21]$}
	\label{fig2}
\end{figure}

The proposed panel data model is applied to estimate the Apple Inc. (AAPL) Call Option price function  using the real market data from August 2019. The call options with the initial maturity of 32 days were observed on a daily basis (except weekends and holidays) until they reached their expiry date (September 6th, 2019).  All together, there are 36 different strikes and the corresponding intrinsic values are fully observed for 21 consecutive day (see Figure \ref{fig1_alt}). The second degree truncated spline basis is used and the overall time-dependent option price function is estimated all-at-once using the minimization formulation in \eqref{surface-estimation} together with the linear constraints in \eqref{constraints}.  The Mosek 9 solver and the R software (Team Development Core, 2019) with the Mosek to R interface (library \texttt{Rmosek}) is used to obtain the solution presented in Figure \ref{fig2}.

It is obvious, that at each time point $t \in \{1, \dots, T\}$, the estimated price function obeys all qualitative properties required for the arbitrage-free market:  the price function is continuous, non-increasing, and convex in the strikes.  The sparsity structure is also obvious as the option price function behaves same in some segments over time: the price function slightly increases after the time $t = 3$ reaching its maximum in between $t = 11$ and $t = 13$. Finally, it again drops down and it stays stable until the end of the follow-up period, $t = 21$. Thus, four changes in the price function are estimated over the follow-up period of 21 days with the time segments  $[1,3]$, $[4,6]$, $[7,10]$, $[11,13]$, and $[14,21]$. 

Let us just briefly mention that due to the spline basis used in \eqref{basis} the panels in Figure \ref{fig2} can be immediately interpolated for any available strike in between two panels. The overall time dependent price surface again satisfies the qualitative restrictions required for the arbitrage-free market. This nicely corresponds with the financial theory and practical applications where the price function should be defined for any arbitrary strike $x \in \mathcal{D}$. The model can be also easily generalized for high frequency data when the intrinsic values are observed on a more dense grind than just on a daily basis. 

Finally, a small simulation was conducted to demonstrate the robustness of the proposed method. The finite sample performace is also compared with the method from \cite{qian_su2}. The follow-up period is $T = 10$ and the sample size $n \in \{20, 100, 200\}$. There is either one changepoint (two time segments) located in the middle of the follow-up period or there are 4 changepoints (5 time segments) equidistantly spaced within the follow-up period. For brevity, the independent error terms are only considered but two error distributions are used: the standard normal distribution and the Cauchy distribution. 

\begin{table}[!t]\footnotesize
	\begin{center}
		\scalebox{0.95}{\input{table_p2.tex}}	
	\end{center}
	\caption{\footnotesize The finite sample comparison of the method from \cite{qian_su2} and the proposed quantile LASSO approach. Two goodness-of-fit quantities are used: the median (MED) of $(Y_{t i} - \widehat{Y}_{t i})$ and the mean absolute difference (MAD) between the true parameter vector $\boldsymbol{\beta}_t^*$ and the corresponding empirical estimate $\widehat{\boldsymbol{\beta}}_t^*$.  The \textit{Recovery} column is given in terms of two values: the proportion of truly discovered  changepoints (value $1$ stands for all true changes being discovered) and the proportion between the number of estimated changepoints and the true changepoints (value $1$ stands for the situation where the number of estimated changes equals the number of true changes). An ideal situation is $1.00/1.00$ which means that all true changes are discovered with no other detections in addition. The results are averaged over 1000 Monte Carlo simulation runs.}
	\label{tab1}
\end{table} 

The performance of  both methods is assessed by using the median of $(Y_{t i} - \widehat{Y}_{ti})$ quantities denoted as MED, the mean absolute difference between the true parameter vectors and their corresponding estimates (denoted as MAD). The changepoint detection is assessed in terms of two quantities: the proportion of truly discovered changepoints (value $1$ stands for all true changes being discovered) and the proportion between the number of estimated changepoints and the true changepoints (value $1$ stands for the situation where the number of estimated changes equals the number of true changes). An ideal situation is $1.00/1.00$ which means that all true changepoints are discovered with no other detections in addition. The results summarized in Table \ref{tab1} are averaged over 1000 Monte Carlo simulation runs.

\section{Conclusion}
\label{conclussion}

The option price function and the implied volatility surface are both fundamental tools for the empirical econometrics, the financial derivatives markets in particular. A new method, based on the panel data structures, is proposed to estimate the time dependent option price function which can be later used to interpolate the implied volatilities. 

 The idea to avoid some standard multistage techniques or nonparametric (semiparametric) smoothing which usually perform slowly and, moreover, additional tuning parameters are required to be specified. Instead, the sparsity principle and the LASSO-type penalty are used to estimate the  option price function which may develop over time. The final model  complies with the arbitrage-free conditions required by the financial theory.

The main advantage of the proposed  method is that it does not apriori assume the arbitrage-free input data. The estimated option price function, which  satisfies the arbitrage-free conditions, is obtained automatically in a~straightforward way by using the estimation procedure together with some well defined linear constraints. These constraints are used to enforce the arbitrage-free scenario. This is crucial for the price computation because the price estimates violating the natural market conditions could have hazardous consequences. 

The estimated option price function can be later used to derive the  implied volatility function (see \cite{Fengler05} or \cite{Kahale04}).  Alternatively, the proposed method can be also directly used to estimate the volatility function using the market observed daily volatilities. Unrestricted market scenarios can be obtained immediately, however, the arbitrage-free scenarios are slightly more complicated and  nonlinear constraints must be used to obtain the final implied volatility function which is compliant with the arbitrage-free conditions derived form the financial theory. 

The proposed quantile LASSO method for the panel data structures serves as an innovative and pioneering approach for the option pricing problem and the following implied volatility estimation. If properly defined, it can effectively handle the estimation under the arbitrage-free criteria which are automatically fulfilled.

\bigskip

\noindent
{\small\textbf{ Acknowledgement: } This work was partially supported by the Czech Science Foundation project GA\v{C}R No. 18-00522Y.}



\appendix
\section{Proof of Theorem \ref{thm1}}\footnotesize
Let us define the quantile proces 
$$G_{n}(\mathbb{B}) \equiv G_{n}(\eb_1, \dots, \eb_T) = \sum_{t = 1}^T \sum_{i = 1}^n \rho_{\tau} (Y_{t i} - \boldsymbol{x}_i^\top \eb_t),$$ 
for $\mathbb{B} = (\eb_1, \dots, \eb_T) \in \mathbb{R}^{p \times T}$ and let $\widetilde{\mathbb{B}}_n \equiv (\widetilde{\eb}_1, \dots, \widetilde{\eb}_T) = \arg \min G_{n}(\eb_1, \dots, \eb_T)$. Firstly, we show that for any $\epsilon > 0$ there is a constant $C_\epsilon > 0$ such that
 \begin{equation*}
\label{eq20}
\PP \left[ \inf_{\textbf{u}_1, \dots, \textbf{u}_T \in \R^{p}, \| \textbf{u}_t\|_1=1 } G_n \big( \mathbb{B}^* + C_\epsilon \delta_n (\textbf{u}_1, \dots, \textbf{u}_T) \big) > G_n(\mathbb{B}^*)   \right] \geq 1 - \epsilon,
\end{equation*}
for $n \to \infty$, where $\mathbb{B}^* = (\eb_1^*, \dots, \eb_T^*)$ are the true parameter vectors and $\delta_n > 0$ is the sequence defined in Assumption (A3).  Equivalently we can write, for any $c > 0$, that 
\begin{align}
 G_n \big( \mathbb{B}^* +  c \delta_n (\textbf{u}_1, \dots, \textbf{u}_T) \big) - G_n(\mathbb{B}^*) & = \E \big[G_n \big( \mathbb{B}^* + c  \delta_n (\textbf{u}_1, \dots, \textbf{u}_T) \big) - G_n(\mathbb{B}^*)  \big] \label{term1}\\
 & \hskip2cm + c \delta_n \sum_{t = 1}^T \sum_{i = 1}^n Z_{ti} \boldsymbol{x}_{i}^\top \textbf{u}_t \label{term2}\\
 & \hskip2cm + \sum_{t = 1}^T \sum^n_{i=1} (W_{ti} -\E[W_{ti}]),\label{term3}
\end{align}
where the random variables $Z_{ti}$ and $W_{ti}$ are defined as $Z_{t i } = (1 - \tau)\mathbb{I}_{\{\varepsilon_{t i} < 0\}} - \tau \mathbb{I}_{\{ \varepsilon_{t i} \geq 0 \}}$ and $W_{t i} = \rho_{\tau}(\varepsilon_{t i} - c \delta_n \boldsymbol{x}_i^\top \textbf{u}_t) - \rho_{\tau}(\varepsilon_{t i}) - c \delta_n Z_{t i}\boldsymbol{x}_{i}^\top \textbf{u}_t$. The random variables $Z_{t i}$ and $W_{t i}$ are independent with respect to $i \in \{1, \dots, n\}$ for some fixed $t \in \{1, \dots, T\}$ but they might be dependent with respect to $t = 1, \dots, T$, where $T \in \mathbb{N}$ is fixed. We need to study the three terms in \eqref{term1}--\eqref{term3}. 

Using the facts that $\rho_{\tau}(x - y) - \rho_{\tau}(x) = y(\mathbb{I}_{\{ x < 0 \}} - \tau) + \int_{0}^y (\mathbb{I}_{\{ x \leq v \}} - \mathbb{I}_{\{ x \leq 0 \}}) \mbox{d}v$, for any $x,y \in \mathbb{R}$, $\E Z_{t i} = 0$ for any $i = 1, \dots, n$ and $t = 1, \dots, T$, the distributional properties from Assumption (A1), and the Taylor expansion, we get
\begin{align*}
\E \big[G_n \big( \mathbb{B}^* + c  \delta_n (\textbf{u}_1, \dots, \textbf{u}_T) \big) - G_n(\mathbb{B}^*)  \big] & = \frac{f(0)}{2}c^2 \delta_n^2 \sum_{t = 1}^T\sum_{i = 1}^n (\boldsymbol{x}_i^\top \textbf{u}_t)^2\\
& + o\left( \delta_n^2 \sum_{t = 1}^T \sum_{i = 1}^n \textbf{u}_t^\top (\boldsymbol{x}_i \boldsymbol{x}_i^\top) \textbf{u}_t\right).
\end{align*}
Thus, due to Assumptions (A1) and (A3), using also the fact that $T \in \mathbb{N}$ is fixed, it also holds that
\begin{align*}
\frac{1}{n}\E \big[G_n \big( \mathbb{B}^* + c  \delta_n (\textbf{u}_1, \dots, \textbf{u}_T) \big) - G_n(\mathbb{B}^*)  \big] & = C f(0) \frac{\delta_n^2}{n}  \sum_{t = 1}^T \sum_{i = 1}^n (\boldsymbol{x}_i^\top \textbf{u}_t)^2 \left(1 + o_P(1)\right) > 0,
\end{align*}
where $C > 0$ is some positive constant.

For the term in \eqref{term2}, we can use the Central Limit Theorem for the independent random variables $\{\sum_{t = 1}^T Z_{t i}\boldsymbol{x}_i^\top \textbf{u}_t\}_{1 \le i \leq n}$ and, again, the fact that $T \in \mathbb{N}$ is fixed, to obtain $$
c \delta_n \sum_{i = 1}^n \sum_{t = 1}^T  Z_{ti} \boldsymbol{x}_{i}^\top \textbf{u}_t = O_P(n^{1/2} \delta_n).
$$
Similarly, for the  term in \eqref{term3}, we have 
\begin{align}
\E \left[ \sum^n_{i=1} \sum_{t = 1}^T(W_{ti} -\E[W_{ti}])\right]^2 & = \sum_{i = 1}^n \E \left[ \sum_{t = 1}^T (W_{t i} - \E [W_{t i}]) \right]^2 \nonumber\\
& \leq \sum_{i = 1}^{n} \sum_{t = 1}^T \E[W_{t i}]^2 + 2 \sum_{i = 1}^{n} \left\{ \sum_{t \neq l} \E [W_{t i} W_{l i}] \right\}, \label{term4}
\end{align}
where the first term in \eqref{term4} can be further bounded by using the moment properties of $W_{t i}$, the distributional assumptions in (A1), and the Taylor expansion for $F(c \delta_n |\boldsymbol{x}_i^\top \textbf{u}_t|)$ and $F(- c \delta_n |\boldsymbol{x}_i^\top \textbf{u}_t|)$ to get 
$$
\sum_{i = 1}^{n} \sum_{t = 1}^T \E[W_{t i}]^2  \leq C \delta_n^2 \sum_{i = 1}^n\sum_{t = 1}^{T} \textbf{u}_t^\top (\boldsymbol{x}_i\boldsymbol{x}_i^\top) \textbf{u}_t \E \big[ \mathbb{I}_{\{ |\varepsilon_{t i}| \leq c \delta_n |\boldsymbol{x}_i^\top \textbf{u}_t| \}} \big] \leq C \delta_n^3 \sum_{i = 1}^n \sum_{t = 1}^T (\boldsymbol{x}_i^\top \textbf{u}_t)^2,
$$
for some positive constant $C > 0$.  Slightly more computational effort is needed to apply the same bound also for the second term in \eqref{term4}. The Davydov's inequality for strictly stationary processes can be used together with the fact that for any $i = 1, \dots, n$ and $t = 1, \dots, T$ it holds that $\E|\mathbb{I}_{\{\varepsilon_{t i} < 0\}}|^{2 + \chi} = \E[\mathbb{I}_{\{\varepsilon_{t i} < 0\}}] = P[\varepsilon_{t i} < 0] = \tau$. Using also the fact that $T \in \mathbb{N}$ is fixed,
we obtain 
$$
\E [W_{t i} W_{l i}] \leq C \delta_n^2 [\textbf{u}_t^\top (\boldsymbol{x}_i \boldsymbol{x}_i^\top)\textbf{u}_l]  \cdot \E \left\{\big[ \mathbb{I}_{\{ |\varepsilon_{t i}| \leq c \delta_n |\boldsymbol{x}_i^\top \textbf{u}_t| \}} \big] \cdot \big[ \mathbb{I}_{\{ |\varepsilon_{l i}| \leq c \delta_n |\boldsymbol{x}_i^\top \textbf{u}_l| \}} \big]\right\},
$$
and the Taylor expansion for the two-dimensional marginal distributions $F_{\varepsilon_{t}, \varepsilon_{t + k}}(x, y)$, for $k = 1, \dots, T - 1$, and Assumption (A1) are needed to get that
$$
\sum_{i = 1}^{n} \sum_{t = 1}^T \E [W_{t i} W_{l i}] = O(n\delta_n^3).
$$

Finally, using the fact that $\frac{1}{n} \sum_{i = 1}^n \textbf{u}_t^\top (\boldsymbol{x}_i\boldsymbol{x}_i^\top) \textbf{u}_l$ is bounded by Assumption (A2) for any $\textbf{u}_t, \textbf{u}_l \in \mathbb{R}^p$, such that $\|\textbf{u}_t\| = \|\textbf{u}_l\| = 1$, and the fact that $T \in \mathbb{T}$ is fixed, we obtain
$$
G_n \big( \mathbb{B}^* +  c \delta_n (\textbf{u}_1, \dots, \textbf{u}_T) \big) - G_n(\mathbb{B}^*) = C n \delta_n^2 \Big( \frac{1}{n}\sum_{i = 1}^n \big[ \sum_{t = 1}^T \textbf{u}_t (\boldsymbol{x}_i \boldsymbol{x}_i^\top) \textbf{u}_t \big] \Big) \cdot (1 + o_P(1)) > 0,
$$
which also implies that $\|\widetilde{\eb}_t - \eb_t^*\|_1 = O_P(\delta_n)$ for any $t = 1, \dots, T$.
To finish the proof we will consider another quantile process 
$$
G_{n}^\star(\mathbb{B}) \equiv G_n(\mathbb{B}) + n \lambda_n \sum_{t = 2}^T \| \eb_t - \eb_{t - 1} \|_2
$$
where, again, $\mathbb{B} = (\eb_1, \dots, \eb_T) \in \mathbb{R}^{p \times T}$. It is easy to see that 
\begin{align*}
G_{n}^\star(\mathbb{B}^* & + c\delta_n (\textbf{u}_1, \dots, \textbf{u}_T)) - G_{n}^\star(\mathbb{B}^*)
 = G_n(\mathbb{B}^* + c\delta_n (\textbf{u}_1, \dots, \textbf{u}_T)) - G_n(\mathbb{B}^*) \\
& \hskip1cm + n \lambda_n \sum_{t = 2}^T \Big[ \| (\eb_t^* + c \delta_n \textbf{u}_t) - (\eb_{t - 1}^* + c \delta_n \textbf{u}_{t - 1}) \|_2  - \|\eb_t^* - \eb_{t - 1}^*\|_2 \Big],
\end{align*}
again for some positive constant $c > 0$. Moreover, from the previous, we already have that 
$$
G_n(\mathbb{B}^* + c\delta_n (\textbf{u}_1, \dots, \textbf{u}_T)) - G_n(\mathbb{B}^*) \geq C n \delta_n^2 > 0,
$$
therefore, using also the triangular inequality, and the fact that $\|\textbf{x} - \textbf{y}\|_2 \leq \|\textbf{x} - \textbf{y}\|_1 \leq \|\textbf{x}\|_1 + \|\textbf{y}\|_1$, for any $\textbf{x}, \textbf{y} \in \mathbb{R}^p$, we have
\begin{align*}
n \lambda_n \sum_{t = 2}^T \Big[ \| (\eb_t^* + c \delta_n \textbf{u}_t) & - (\eb_{t - 1}^* + c \delta_n \textbf{u}_{t - 1}) \|_2  - \|\eb_t^* - \eb_{t - 1}^*\|_2 \Big] \\
& \geq n \lambda_n \sum_{t \in \mathcal{A}^*} \Big[ \| (\eb_t^* + c \delta_n \textbf{u}_t)  - (\eb_{t - 1}^* + c \delta_n \textbf{u}_{t - 1}) \|_2  - \|\eb_t^* - \eb_{t - 1}^*\|_2 \Big]\\
& \geq - c n \lambda_n \delta_n \sum_{t \in \mathcal{A}^*} \|\textbf{u}_t - \textbf{u}_{t -1}\|_2 = - C n \lambda_n \delta_n,
\end{align*}
for some positive constant $C > 0$, where $\mathcal{A}^* = \{t \in \{2, \dots, T\};~\eb_t^* \neq \eb_{t - 1}^*\}$. Using now Assumption (A3)  we conclude that 
$$
G_{n}^\star(\mathbb{B}^*  + c\delta_n (\textbf{u}_1, \dots, \textbf{u}_T)) > G_{n}^\star(\mathbb{B}^*),
$$
which holds with probability converging to one, as $n \to \infty$. This also implies, that 
$$
\|\widehat{\eb}_t - \eb_{t}^*\|_1 = O_{P}(\delta_n), 
$$
for any $t \in \{1, \dots, T\}$, which completes the proof of Theorem \ref{thm1}. \qed

\end{document}

%% file: table_p2.tex
\begin{tabular}{ccc|cccccc}
\multirow{2}{*}{$\boldsymbol{\mathcal{D}}$} & & \scalebox{0.8}{$\boldsymbol{T = 10}$} & 
     \multicolumn{3}{c}{\textbf{Model with 2 phases}}  & \multicolumn{3}{c}{\textbf{Model with 5 phases}}\\
 ~ & ~ & $\boldsymbol{n}$  & 
     \multicolumn{1}{c}{\scalebox{0.8}{MED}} & \multicolumn{1}{c}{\scalebox{0.8}{MAD}} & \multicolumn{1}{c}{\scalebox{0.8}{Recovery}} &
     \multicolumn{1}{c}{\scalebox{0.8}{MED}} & \multicolumn{1}{c}{\scalebox{0.8}{MAD}} & \multicolumn{1}{c}{\scalebox{0.8}{Recovery}} \\\hline\hline
\multicolumn{3}{c}{~} & \multicolumn{6}{c}{~}\\[-0.3cm]

$\boldsymbol{N}$ & & \textbf{20} & -0.03 & 0.17 & 1.00$/$9.87&0.02 & 1.00 & 0.65$/$1.11\\
 & \scalebox{0.8}{Model from \cite{qian_su2}} & \textbf{100} & 0.00 & 0.05 & 1.00$/$62.67&0.00 & 0.08 & 1.00$/$6.72 \\
 & & \textbf{200} & 0.00 & 0.02 & 1.00$/$87.79&-0.07 & 0.07 & 1.00$/$11.03 \\
 
\multicolumn{3}{c}{~} & \multicolumn{6}{c}{~}\\[-0.2cm]
\rowcolor{Gray} \multicolumn{3}{c}{~} & \multicolumn{6}{c}{~}\\[-0.2cm]

\rowcolor{Gray} &  & \textbf{20} & -0.04 & 0.18 & 1.00$/$9.92&0.08 & 1.01 & 0.57$/$1.02\\
\rowcolor{Gray} & \scalebox{0.8}{Model from \eqref{surface-estimation}} & \textbf{100} & -5.79 & 2.05 & 1.00$/$99.00&-4.86 & 0.91 & 1.00$/$10.08 \\
\rowcolor{Gray} & & \textbf{200} & -0.01 & 0.02 & 1.00$/$29.92&-0.02 & 0.03 & 1.00$/$4.73 \\

\rowcolor{Gray} \multicolumn{3}{c}{~} & \multicolumn{6}{c}{~}\\[-0.2cm]
\rowcolor{white} \multicolumn{3}{c}{~} & \multicolumn{6}{c}{~}\\[-0.2cm]

$\boldsymbol{C}$ & & \textbf{20} & 0.00 & 3.24 & 0.73$/$10.85&0.01 & 3.54 & 0.65$/$1.21\\
 & \scalebox{0.8}{Model from \cite{qian_su2}} & \textbf{100} & -0.08 & 3.36 & 0.86$/$71.60&-0.22 & 3.12 & 0.82$/$7.64 \\
 & & \textbf{200} & -0.03 & 2.29 & 0.79$/$128&-0.12 & 2.30 & 0.80$/$14.52 \\
 
\multicolumn{3}{c}{~} & \multicolumn{6}{c}{~}\\[-0.2cm]
\rowcolor{Gray} \multicolumn{3}{c}{~} & \multicolumn{6}{c}{~}\\[-0.2cm]

\rowcolor{Gray} &  & \textbf{20} & -0.02 & 0.55 & 0.76$/$9.84&0.05 & 1.17 & 0.56$/$1.09\\
\rowcolor{Gray} & \scalebox{0.8}{Model from \eqref{surface-estimation}} & \textbf{100} & -8.20 & 2.03 & 1.00$/$97.88&-7.85 & 1.62 & 0.98$/$10.38 \\
\rowcolor{Gray} & & \textbf{200} & -2.66 & 0.39 & 1.00$/$109&-1.27 & 0.22 & 1.00$/$6.18\\

\rowcolor{Gray}\multicolumn{3}{c}{~} & \multicolumn{6}{c}{~}\\[-0.2cm]\hline\hline

     \rowcolor{white}\multicolumn{3}{c}{~} & \multicolumn{6}{c}{~}\\[-0.4cm]\hline

     \end{tabular}